\documentstyle[12pt]{article}
\newtheorem{defin}{Definition}
\newtheorem{prop}{Proposition}

\newtheorem{th}{Theorem}
\newtheorem{lemma}{Lemma}
\newtheorem{cons}{Corollary}

\newfont{\sdbl}{msbm9}
\newfont{\dbl}{msbm10 at 12pt}
\newcommand{\eqdef}{\stackrel{\rm def}{=}}
\newcommand{\proof}{{\bf Proof\ }}
\newcommand{\eq}{=}
\newcommand{\oo}{{\cal O}}
\newcommand{\ff}{{\cal F}}

\newcommand{\Lim}{\mathop {\rm lim}}

\newcommand{\Proj}{\mathop {\rm Proj}}
\newcommand{\Frac}{\mathop {\rm Frac}}
\newcommand{\dm}{\mathop {\rm dim}}
\newcommand{\trdeg}{\mathop {\rm trdeg}}

\newcommand{\dz}{{\mbox{\dbl Z}}}

\newcommand{\dn}{{\mbox{\dbl N}}}
\newcommand{\sdn}{{\mbox{\sdbl N}}}

\begin{document}
\author{D.V. Osipov, A.B. Zheglov\footnote{supported
by DFG-Schwerpunkt ''Globale Methoden in der Komplexen Geometrie'',
by RFBR grant no. 05-01-00455, by grant of Leading Scientific
Schools no. 489.2003.1, by INTAS grant 05-100000-8118; besides, the
first author was supported by RFBR grant no. 04-01-00702 and by
grant of Russian Science Support Foundation. } }

\title{On some questions related to the Krichever correspondence}
\date{}
\maketitle

\begin{abstract}
We investigate various new properties and examples of
two-dimensional and one-dimensional Krichever correspondence.
\end{abstract}

\section{Introduction}

The aim of this paper is to show several new properties and examples
of the  two-dimensional Krichever correspondence. On the one hand,
these examples  make a general picture of the Krichever
correspondence more clear.  On the   other hand, they should show a
possible way how to generalize the one-dimensional Krichever
correspondence  to the case of dimension 2.

Also we give some explicit examples for the one-dimensional
Krichever correspondence, which explain some unclearness in the
paper~\cite{Mu}, and give some new KP-equations which appear as
corollaries from generalized  KP-hierarchy on two-dimensional local
skew fields. These skew fields were first considered in~\cite{P3}
and then classified in~\cite{Zh1}.

First, we recall that the Krichever map for algebraic curves is the
following construction. We start from a data  $(C, p, \ff, e_p, t)$,
where $C$ is an irreducible algebraic curve, $p$ is a smooth
$k$-point on $C$, $t$ is a local formal parameter at $p$, $\ff$ is a
torsion free coherent sheaf of rank~$r$ on $C$, $e_p$ is a local
formal trivialization of the sheaf $\ff$ at $p$. We consider the
embedding from $H^0 (C \setminus p, \ff)$ into $K_p \otimes \ff $.
This embedding gives us the map from our data to Fredholm subspaces
of the space $ k((t))^{\oplus r}$ by  identifying $K_p \otimes \ff
\simeq K_p^{\oplus r}$ under fixed $e_p$ and $t$. This map is called
the Krichever map, see~\cite{M}, \cite{SW}, \cite{O1}.

Now we recall  the Krichever map for algebraic surfaces
 (see~\cite{P4}, \cite{Par}).  We  start from a data $(X, C, p,
\ff, e_p, t, u )$, where $X$ is a projective irreducible normal
algebraic surface, $C$ is an ample Cartier divisor  on $X$, $p$ is a
smooth  $k$-point on $X$ and $C$, $\ff$ is a rank $r$ vector bundle
on $X$, $e_p$ is a local formal trivialization of $\ff$ at $p$ on
$X$, $u$ and $t$ are local parameters at $p$ such that $u=0$ is a
local equation of $C$ on $X$ near $p$. By such data a $k$-subspace
$W$ in the sum of two-dimensional local fields $k((t))((u))^{\oplus
r}$ is canonically constructed.

Shortly, this construction consists of  the following steps. By the
data $p \in C \subset X$ we can canonically build the
two-dimensional local field $K_{p,C}$, which is isomorphic to
$k((t))((u))$ after fixing local parameters $u$ and $t$ (see
\cite{P1}, \cite{PF}, \cite{H}, \cite{O2}). We consider the subfield
$K_C\subset K_{p,C}$, which is the completion of the function field
on $X$ with respect to the valuation given by the curve $C$, and the
subring $B_p = \Lim\limits_{\longrightarrow n}
\frac{\hat{\oo}_p}{u^n}$, where $\hat{\oo}_p$ is the complete local
ring at $p$  on $X$. ($B_p$ does not depend on the choice of $u$).
Now we have
$$
W =  \bigcap_{q \in C, \; q \ne p} (B_q \otimes \ff) \cap (K_C \otimes \ff)
$$
And from the embedding of the space $K_C \otimes \ff$ to the space
$K_{p,C} \otimes \ff $ we have the embedding of the space $W$ to the
space $k((t))((u))^{\oplus r}$ by chosen $e_p$, $t$, and $u$.

Note also that the Krichever map was generalized to  algebraic
varieties of higher dimensions in the work~\cite{O3} by
higher-dimensional local fields.

We are grateful to A.N. Parshin for the stimulating remarks and
discussions. We are grateful also to H. Kurke, and to the Humboldt
University of Berlin for the nice atmosphere, where the part of this
work was done.

\section{Some explicit examples of the two-dimensional Krichever correspondence}
In this section we shall compute some examples of the two-dimensional
Krichever correspondence on $\bf P^2$.

In~\cite{P4} was given an example when $X = \bf P^2$ with
homogeneous coordinates $(x:y:z)$, the curve $C$ is given by $z=0$,
the point $p$ is $(1:0:0)$, the bundle $\ff = \oo_{X}$, the local
parameters $t= y/x$, $u = z/x$. Then the subspace $W =
k[t^{-1}]((t^{-1} u))$.

For our examples we need the following lemma.

\begin{lemma}  \label{le}
Suppose that there is an ample effective  divisor $D$ on $X$ such
that $D \cdot C = m p$ for some integer $m$, and the affine curve $C
\setminus p$ is given by the equation $f = 0$ on the affine surface
$X \setminus D$ for some $f \in k[X \setminus D]$.
 Choose any $y_1, \ldots, y_l \in  k[X \setminus D]$
such that their images in $k[C \setminus p ]$ generate this ring.
(Such $y_i$ always exist). We have that canonically  $y_1, \ldots,
y_l, f \in K_{p, C}$. Then we have
$$W = k[y_1, \ldots, y_l ] ((f))  \subset K_{p,C}$$.
\end{lemma}
\proof follows very easily from the definition of subspace $W$ and
facts that $K_C k(C)((f))$ and  $B_q = \Lim\limits_{\longrightarrow
n} \frac{\hat{\oo}_q}{f^n}$ for any point $q \in C$, $q \ne p$.

\vspace{0,5cm}

Now we shall apply this lemma to the case $X= \bf P^2$.

Consider the case of quadric curve.
\begin{prop} \label{pe1}
Consider the surface $X = \bf P^2$ with the homogeneous coordinates
$(x : y : z)$. Let the curve $C$ be given by the equation $y^2 + z^2
- 2xz = 0$, and the point $p= (1:0:0)$. Let local parameters $t =
y/x$, $u = (y/x)^2 + (z/x)^2 - 2 z/x$. Then for the sheaf $\ff =
\oo_X$ we have
$$
W = k[\alpha] (( \frac{\alpha^2}{t^2} u)) \mbox{,}
$$
where
$ \alpha = \frac{t (1 + (1 +u - t^2)^{1/2})}{t^2 - u}$.
\end{prop}
\proof. To apply lemma~\ref{le}, we consider the curve $D$ which is
the tangent line to the curve $C$ at the point $p$. So $D$ is given
by the equation $z = 0$, and $D \cdot C = 2p$. Now in the affine
domain $z \ne 0$ we have an equation $f = (y/z)^2 +1 - 2 x/z$. Note
that $f = u x^2/ z^2$, and $k[C \setminus p] = k[y/z]$. Denote
$\alpha = y/z$. Then from lemma~\ref{le} we obtain that $W =
k[\alpha] (( u \alpha^2 / t^2 ))$. Let us calculate $\alpha$. From
the expression for $u$ we have that $ u = t^2 + (t/\alpha)^2 - - 2 t
/ \alpha $. Hence $t / \alpha = 1 + (1+u - t^2)^{1/2}$ or $t /
\alpha = 1 - (1+u - t^2)^{1/2}$. But $t / \alpha  = z/x$ has zero of
order~$2$ at the point $p$ after the restriction to the curve $C$.
Therefore for $t / \alpha$ we have the second variant, and hence we
obtain the formula from the proposition. The proposition is proved.

\vspace{0.5cm}

Note that for any quadric on $\bf P^2$ we can take the tangent to
the point and apply lemma~\ref{le}.

Now we consider the case of any curve $C$  of degree $m$ on $\bf
P^2$. Suppose that there exists a point $p \in C$  such that the
tangent line $D$ to the curve $C$ at $p$ has the property $D \cdot C
= mp$.

For example, if $C$ is a cubic curve, then such point $p$ always
exists. (But for the case of a quartic curve this fact is not true
anymore, for a general quartic curve such point does not exist.)

Choose homogeneous coordinates $(x:y:z)$ on $\bf P^2$ such that $D$
is given by $z = 0$, and $p$ is $(1:0:0)$.  Then we have that the
affine curve $C \setminus p$ on the affine surface $X \setminus D$
is given by the equation $f (x/z , y,z) = 0$. Therefore we can
choose the local parameters $u = (z/x)^m f$ and $t = y/x$. Also the
images of elements $x/z$ and $y/z$ generate the ring $k [C \setminus
p]$. Hence by lemma~\ref{le}  for the sheaf $\ff = \oo_X$ we obtain
\begin{equation} \label{ex1}
W = k[ \beta, \beta t] (( {\beta}^m u ))  \mbox{,}
\end{equation}
where $\beta = x/z$ is defined from the equation $f(\beta, \beta t)
= 0$ and the condition that $\beta$  has a pole of order $m$ after
the restriction to the curve $C$.

{\em Define} the support of a $k$-subspace $W$ from $k((t))((u))$ as
the closed $k$-subspace of $k((t))((u))$ generated by all monomials which
are the lowest degrees of elements from  $W$.

Note that the similar subspaces play important role
 for the case of Krichever correspondence for algebraic curves.
For example, the "holes" between monomials of negative degrees are
Weierstrass gaps, their quantity is the genus of the curve. Also in
this case it is possible to reconstruct the cohomology groups of
original sheaf on the curve from the support of a subspace arising
from a geometrical data.

If we are interesting in  $k$-subspaces from $k((t))((u))$
 which are images of the Krichever correspondence for algebraic surfaces,
it makes sense to investigate the possible supports for such
subspaces. But note at once that not all the properties will hold
for  the case of algebraic surfaces. So, from the support of the
subspace we can not reconstruct the  cohomology groups of the
original sheaf on a surface. (For example, we do not know   anything
about the cohomology group $H^1$ from the support of the subspace.)

Now from explicit calculations above we have some examples of subspaces-supports
which appear from the Krichever correspondence for algebraic surfaces.
\begin{prop} \label{pe2}
Let the conditions be as in formula~\ref{ex1}. Then the support of
such subspace $W$ is the closed subspace generated by all monomials
$t^i u^j$ with $i \le -m^2j $ except for the monomials
$t^{-\alpha_1-m^2l} u^{l}, \ldots, t^{-\alpha_g-m^2l} u^{l}$, where
$i$, $j$  and $l$ are any integers,  and  natural numbers $\alpha_1,
\ldots , \alpha_g$ are the Weierstrass gaps for the curve from our
data.
\end{prop}
\proof follows by the direct application of formula~\ref{ex1}. And
$\alpha_1, \ldots, \alpha_g $ are natural numbers which can not be
obtained by linear combinations of integers $m$ and $m-1$ with
integer nonnegative coefficients ($g$ is the genus of the curve).

\section{Stabilizer ring of a k-subspace of a two-dimensional local field.}

This section is one of the first steps to describe explicitly the images of the
two-dimensional Krichever map, that is, to determine $k$-subspaces
of $k((t))((u))$ which are in the image of the two-dimensional Krichever map.

Recall the case of the Krichever correspondence for the curves
(see~\cite{M}). For any Fredholm $k$-subspace of $k((t))$, that is,
for the $k$-subspace $W$ such that  $\dm_k W \cap \oo < \infty$, and
$\dm_k k((t))/ (W + \oo) < \infty$, where $\oo = k[[t]]$, we
consider the ring
$$A \eqdef \{a \in k((t)) : a \cdot W \subset W\} \mbox{.} $$
Then the ring $A$ is finitely generated over the field $k$, and
either $A = k$, or the ring $A$ has the Krull dimension equal to
$1$.

In the second case the subspace $W$ is in the image of the Krichever
map, that is, the curve $C = \Proj (\mathop{\oplus}\limits_n A \cap
t^{-n}\oo )$, the sheaf $\ff = \Proj (\mathop{\oplus}\limits_n A
\cap t^{n} \oo )$, the point $p$ appears from the given filtration
$A \cap t^{n} \oo$ of the ring $A$, $t$ is a formal local parameter
at $p$, and $1 \in k((t))$ is a formal trivialization of $\ff$.
Moreover, the original $k$-subspaces $A$ and $W$ appear as the
images of the Krichever map of the constructed geometrical data $(C,
p, \oo_C, 1, t)$ and $(C, p, \ff, e_p, t)$ correspondingly. Note
that this description is connected with the description of
commutative subrings in the ring of differential operators in one
variable, that is, with the "Schur pairs", see~\cite{M}.

Recall that the Euler characteristic of a Fredholm subspace
$W\subset k((t))$ is equal to $\chi(W)\eq \dm_k W \cap \oo - \dm_k
k((t))/ (W + \oo )$.

Consider now the case of a two-dimensional local field. Recall that
a two-dimensional local field has, by definition, two valuations
$\nu$ and $\bar\nu$, where  $\nu$ is a discrete valuation of rank 1
on a two-dimensional local field, and $\bar\nu$ is a discrete
valuation of rank 1 on the residue field of a two-dimensional local
field. Let the space $V \eq k((t))((u))$, and define $k$-subspaces
$\oo_1 \eq k((t))[[u]]$, \ $\oo_2 \eq k[[t]]((u))$. For any integer
$n$, for any $k$-subspace $W \subset V$ let $W(n) = (t^n \oo_1 \cap
W) / (t^{n+1}\oo_1  \cap W)$.

For any commutative ring $B \supset k$ we say $\trdeg_k B= m$ if $m$
is the maximal number of algebraically independent elements of the
ring $B$ over the field $k$.

\begin{th}
\label{main}

Let $W$ be a $k$-subspace of the space $V$ such that for any integer
$n$ the space $W(n)$ is a Fredholm subspace in a one-dimensional
local field, and $\chi (W(n)) = a + bn$, where  $b < 0$. Let the
ring $A$ be a $k$-subring of the space $V$ such that $A \supset k$,
\ $A \cdot W \subset W$. Then

a) for any element $a\in A$ we have $\bar\nu
(a)\le b\nu (a)$.

b) $\trdeg_k \Frac (A \cap \oo_2) \le 2 $, and the field $\Frac (A
 \cap \oo_2)$ is finitely generated over the ground field $k$.
(For any commutative ring without zero divisors by $\Frac B$ we
denote its   field of fractions.)
\end{th}

{\bf Proof.} a) We assume the converse. Then there exists an element
$x\in A$ such that $\bar\nu (x)> b\nu (x)$. We have $x\cdot W\subset
W$ and $x\cdot W(0)\subset W(m)$. It is easy to see that $\chi
(x\cdot W(0))\eq \chi (W(0))+\bar\nu (x)$. Now we have
$$
\chi (W(m))\eq a+bm<a+ \bar\nu (x)\eq \chi (W(0))+\bar\nu (x) \eq
\chi (x\cdot W(0))\le \chi (W(m)).
$$
It is a contradiction.

b) Instead of this case we will prove  more general result:

\begin{lemma}
Let $B$ be a subring in a two dimensional local field $V \eq k((t))((u))$ such
that $k \subset B$ and the following condition holds:

for every $a\in B$ we have $0  \le \bar\nu (a) \le - \pi \nu (a)$, \
$\pi
> 0$.

Then $\trdeg_k \Frac B \le 2 $, and the field $\Frac B$ is finitely
generated over the ground field $k$.
\end{lemma}

{\bf Proof.} Consider the subspace $\bigtriangleup (N)\eq \{a\in B |
\nu (a)> - N\}$, where the integer $N > 0$. Note that this subspace
has a finite dimension over the field $k$, and this dimension is not
greater than $\pi N^2$.

 Indeed, if this dimension is greater than $\pi N^2$, then the ring $B$ must contain an
element $x$ such that $\bar \nu (x) > -\pi \nu (x)$. It is a
contradiction.

Note that $\trdeg_k B< 2 $  if and only if $\bar\nu (a) \eq l \nu
(a)$ for some constant $l$ and all elements $a\in B$. Indeed, if
there are two elements $a,b$ such that the last condition does not
hold, then they must be algebraically independent, because the first
monomials of such elements are algebraically independent. Converse,
if $\bar\nu (a) \eq l \nu (a)$ for all $a \in B$, then the supports
(the monomials of lowest degrees) of all elements from $B$ lie on a
line which goes through the zero point of a coordinate plane, and
therefore one can apply arguments, for example, from \cite[Prop.
3.2.]{M} to obtain $\trdeg_k B< 2$.

We suppose now that $\trdeg_k B \ge 2 $. We take two elements $a,b
\in B$ such that $\bar\nu (a)/ \nu (a) \ne  \bar\nu (b)/ \nu (b)$.
It means that their supports do not belong to a line which goes
through the zero point on a coordinate plane. These elements exist
because $\trdeg_k B \ge 2 $, and because of the arguments from the
previous paragraph. Also, they are algebraically independent
elements. We can estimate the dimension of a subspace generated over
the field $k$ by powers of elements $a,b$ which lie in the subspace
$\bigtriangleup (N)$ for some $N$.
 It is clear that this dimension is greater than
 $$  \frac{1}{2} \cdot (-\frac{N}{\nu(a)} -1 )(-\frac{N}{\nu(b)} -1) =  eN^2+lN+r $$
 for some real numbers
$e,l,r$ and for every sufficiently large number $N$. (To obtain this
estimate we consider a parallelogram constructed by vectors which go
from the zero point of a coordinate plane and end in the points
$a^{[-N /\nu(a)]}$ and $b^{[-N /\nu(a)]}$ correspondingly.)

Now suppose that the field $\Frac B$ has infinite dimension over the
field $\Frac k[a,b]$. Then we take $M$ linearly independent over the
ring $k[a,b]$ elements $a_1,\ldots, a_M$, where the number $M$
satisfies the condition $Me> \pi$. Without loss of generality we can
assume
 that
$\nu (a_1)\le \nu (a_i)$ for all $i$. Since the dimension of a
subspace generated over the field $k$ by powers of elements $a,b$
which lie in the subspace $\bigtriangleup (N+M\nu (a_1))$ is greater
than $e(N+M\nu (a_1))^2+l(N+M\nu (a_1))+r$ for every sufficiently
large number $N$, we obtain that the dimension of a subspace
generated over the field $k$ by elements $a_1,\ldots a_M$ multiplied
by powers of elements $a,b$ which lie in the subspace
$\bigtriangleup (N+M\nu (a_1))$ is greater than $M(e(N+M\nu
(a_1))^2+l(N+M\nu (a_1))+r)$ for every sufficiently large number
$N$.  The last subspace is inside the space $\bigtriangleup (N)$. On
the other hand, since $Me>\pi$,  for every sufficiently large number
$N$ we have
$$M(e(N+M\nu (a_1))^2+l(N+M\nu (a_1))+r)> \pi N^2>
\dm\nolimits_k \bigtriangleup (N) \mbox{.}$$ We obtained a
contradiction. Therefore, in this case the field $\Frac B$ has
finite dimension over the field $\Frac k[a,b]$. Hence, $\trdeg_k
\Frac B = 2 $, and the field $\Frac B$ is finitely generated over
the ground field $k$.

The lemma is proved.

The theorem is proved.

\vspace{0.5cm}

{\bf Remark 1.} As we have seen in the proof of the lemma, $\trdeg_k
B \eq 1 $ if and only if the supports of all elements from the ring
$B$ lie on a
 line which goes through zero point on a coordinate plane.
(Moreover, arguing similar to~\cite[Prop. 3.2.]{M}), it is possible
to prove that the ring $B$ is finitely generated over the field $k$
and has Krull dimension equal $1$.)

  If
$\trdeg_k B\eq 0 $, then $B \eq k$.

\vspace{0.5cm}

{\bf Remark 2.} It is important that in theorem~\ref{main} we
consider the intersection $A \cap \oo_2$. The statement b) of
theorem~\ref{main} is not true for the ring $A$ itself as the
following example shows.

Consider the subspace $W \eq \{ a\in V | \bar\nu (a) \le -\nu (a)
\}$. It is possible to check that the subspace $W$ satisfies the
conditions of theorem~\ref{main}, and that $A \eq W$. But since the
ring $A$ contains the subfield $k((ut^{-1}))$, we have that the
field $\Frac A$ has the infinite transcendency degree over the field
$k$.

\vspace{0.5cm}

{\bf Remark 3.} In theorem~\ref{main} we consider  the ring $A \cap
\oo_2$, because this ring is the image of the affine ring $H^0 (X
\setminus C, \oo_X)$ when the pair $(A, W)$ is the image of the
geometrical data $(X, C, p, \oo_X, 1, t, u)$, \ $(X, C, p, \ff, e_p,
t, u)$ under the Krichever map. And we need exactly this ring for
the  theorem of reconstruction of geometrical data. Namely, the
required surface  $X = Proj (\mathop{\oplus}\limits_n (A \cap \oo_2)
\cap t^n \oo_2)$, the required curve $C$ appears from the filtration
induced by the filtration $t^n \oo_2$ on $A \cap \oo_2$, and the
required point $p$ appears from the one-dimensional Krichever
correspondence on the curve $C$, which is given by $W(0)$
(see~\cite{P4}).

Also the space  $W(n)$ is the image of the sheaf $\ff \mid_C \otimes
N_{C/X}^{\otimes -n}$ under the Krichever map on the curve~$C$. And
by the Riemann-Roch theorem applied to the curve $C$ we obtain the
Euler characteristic
$$\chi (N_{C/X}^{\otimes -n}) = 1 - g(C) + deg
(\ff) - n  \, deg N_{C/X} \mbox{.} $$ Besides $N_{C/X} = \oo (C
\cdot C) $. And from ampleness of the divisor $C$ we have $deg
N_{C/X}
> 0 $. It explains the condition $\chi (W(n)) = a + bn$, \ $b < 0$
from theorem~\ref{main}.

\vspace{0.5cm}

{\bf Remark 4.} If a $k$-subspace $W$ is the image of a geometrical
data  $(X, C, p, \ff, e_p, t, u)$ under the Krichever map, then from
the complex which calculates the cohomology groups (see~\cite{P4},
\cite{Par}) we have
$$ W\cap {\cal O}_1\cap {\cal O}_2 = H^0 (X, \ff)$$
 $$ \frac{W\cap ({\cal O}_1+ {\cal O}_2)}{W\cap {\cal O}_1 + W\cap {\cal
O}_2} = H^1 (X, \ff)$$
 $$ \frac{V}{W+ {\cal O}_1 + {\cal O}_2} = H^2 (X, \ff) \mbox{.}$$
 And we know that the dimensions of these subspaces over the field $k$ are finite
 as  dimensions of cohomology groups.
 It would be interesting to add these
 conditions of finite-dimensionality
 to conditions of theorem~\ref{main} and  to obtain new
 corollaries for the answer, if a $k$-subspace $W \subset V$
 is in the
 image of geometrical data under the Krichever map?

\section{To the question about an embedding of a projective line into an infinite-dimensional
affine space}

In the paper \cite[p. 13]{Mu} was remarked the following curious
fact: there is an embedding of the universal Sato Grassmannian
into an infinite-dimensional affine space. In this section we
clarify, if this embedding is algebraic. To this end we
investigate an example, in which this embedding are restricted to
the projective line in the universal Sato Grassmannian.

First we recall some well-known facts about the Sato Grassmannian,
and also the construction of an embedding map from the
paper~\cite{Mu}.

Let $V \eq k((z))$ be the field of Laurent power series with a
filtration $V(n) \eq z^n k[[z]]$.
 Let $V_1\eq V(0)$.

By $Gr(V)$ we denote the set of  subspaces $W$ in $V$ such that the
complex
$$
W\oplus  V_1\rightarrow V
$$
is a Fredholm one. This set possess the structure of an
infinite-dimensional projective variety, whose connected components
are marked by the value of the Euler characteristic of the complex.
We will work only with "zero" component $Gr_0(V)$.

Now we consider  the skew field $P\eq k((x))((\partial^{-1}))$ of
formal pseudo-differential operators with coefficients from the
field $k((x))$. This skew field is  a left $k((x))$-module of all
expressions $L\eq \sum_{i>\infty }^n a_i\partial^i$, $a_i\in k((x))$
with a multiplication defined according to the Leibnitz rule.

In the same way one can define the ring $E\eq
k[[x]]((\partial^{-1}))$, $E\subset P$. It can be checked that $P$
and $E$ are associative rings (see  details in the paper~\cite{P3}).

There is a decomposition
$$
E\eq E_{+}+E_{-},
$$
where $E_{-}\eq \{ L\in E: L\eq \sum_{n<0}a_n\partial^n \}$, and
$E_{+}$ consists of the operators containing only nonnegative powers
of $\partial$. The elements from $E_{+}\eq :D$ are the usual
differential operators, and the elements from $E_{-}$ are the
Volterra operators.

The map $E\rightarrow E/Ex\eq V$ (we identify the image of
$\partial^{-1}$ with $z$) defines a linear action of the ring $E$ on
$V$ and therefore on $Gr(V)$. The map $E \rightarrow V$ is called
{\em the Sato map}.

Let us introduce the notion of a standard subspace:

\begin{defin}
Let $S$ be a subset in $\dz$, $S:\eq \{\sigma(0), \sigma (-1),
\sigma (-2), \ldots \}$, where $\sigma (-i)\in \dz$, $\sigma (-i)\eq
-i$ for $i\gg 0$.

The subspaces $V^S:\eq \bigoplus_{l \in \sdn \cup 0} k\cdot
z^{\sigma (-l)}\subset V$, $V^S\in Gr_0(V)$ are called the standard
subspaces.

 By $W_0$  we  denote the subspace $V^{S_0}$, where $S_0 \eq \{ -1,-2,\ldots
\}$.
\end{defin}

The following lemma was proved in the paper~\cite{Sa1}.

\begin{lemma}
There is a unique operator $R$ such that $RV^S \eq W_0$. Namely,
$$
R \eq \partial^{-m-1} (x\partial - \sigma(0))(x\partial -(\sigma
(-1)+1))(x\partial -(\sigma (-2)+2))\ldots  (x\partial -(\sigma
(-m)+m)) \mbox{,}
$$
where $m$ is the maximal number such that $\sigma (-m)-m \ne 0$.
\end{lemma}

The notion of a quasiregular operator was introduced in the
paper~\cite{S}.

\begin{defin}
Let ${\cal E}_{K}^{(0)}$, $K\eq k((x))$ be a group of monic (i.e.,
with the leading coefficient equal $1$) operators of degree zero
from the ring $P$.

An operator $W\in {\cal E}_{K}^{(0)}$ is called quasiregular if
there are numbers $m,n\in \dn$ such that $x^mW$ and $W^{-1}x^n$
belong to $E$.

By $\cal R$  we denote the set of all quasiregular operators.
\end{defin}

Now we can introduce a map  $\gamma : {\cal R} \to Gr_0(V)$. Put
$$
\gamma (W)\eq (W^{-1}x^n)W_0, \mbox{\quad} W\in {\cal R}  \mbox{.}
$$
Since $x \cdot W_0 \eq W_0$, this definition is correct, i.e., it
does not depend on $n$. It is clear that $\gamma (W)\in Gr_0(V)$.

The following theorem was proved in the paper~\cite{Sa1}.

\begin{th}
The map $\gamma : {\cal R} \rightarrow Gr_0(V)$ is a bijection.
\end{th}

It was remarked in  the paper \cite[p. 13]{Mu} that ${\cal
E}_{K}^{(0)}$ is an infinite dimensional affine space, and $\cal
R$ is embedded into this space. (It is more correctly to say that
${\cal E}_{K}^{(0)}$ is an inductive limit of infinite-dimensional
affine spaces). It was also remarked that due to the theorem there
is an embedding of the universal Sato Grassmannian into this
affine subspace. To check this assertion we give below an explicit
computation  of this embedding for a projective line which lies in
the  Sato Grassmannian.

Consider the set of Fredholm subspaces
$$
R \eq \{W(\alpha ,\beta )\subset V, W(\alpha ,\beta )\eq
\bigoplus_{l\eq 1}^{\infty} k\cdot z^{-l}\bigoplus k\cdot (\alpha
+\beta z)\},
$$
where $\alpha ,\beta \in k$. It is clear that $R$ is a projective
line in $Gr_0(V)$ with coordinates $(\alpha :\beta )$.

Consider the operators $\tilde{S}(\alpha ,\beta )\eq \alpha +\beta x+\beta
\partial^{-1}$. We have  $\tilde{S}(\alpha ,\beta )\in E$. Note that
$\tilde{S}(\alpha ,\beta )W_0\eq  W(\alpha ,\beta )$ for any pair
$(\alpha ,\beta )$. In the case $\alpha\ne 0$ we define the operator
$S(\alpha ,\beta )\eq \alpha^{-1}\tilde{S}(\alpha ,\beta )(1+\beta
/\alpha x)^{-1}$. Since $(1+\beta /\alpha x)W_0\eq W_0$, we have
$S(\alpha ,\beta )W_0\eq W(\alpha ,\beta )$. Note that $S\in {\cal
E}_{K}^{(0)}\cap E$. In the case $\alpha\eq 0$ it is possible to
represent the operator $\beta^{-1}\tilde{S}(\alpha ,\beta )$ in the
form $H(\alpha ,\beta )^{-1}x$, where $H(\alpha ,\beta )\in {\cal
E}_{K}^{(0)}$, $H(\alpha ,\beta )\eq (1+\partial^{-1}x^{-1})^{-1}\eq
1-x^{-1}\partial^{-1} +\ldots$. Thus, we obtain
$$
\gamma^{-1} (W(\alpha ,\beta ))\eq \left\{ {S(\alpha ,\beta )^{-1}
\eq 1- (\frac{\beta}{\alpha} -\frac{\beta}{\alpha} x
+\frac{\beta^2}{\alpha^2}x^2+\ldots )\partial^{-1} +\ldots
\newline \mbox{\quad if \quad } \alpha\ne 0} \atop {H(\alpha
,\beta )\eq 1-x^{-1}\partial^{-1} +\ldots \mbox{\quad if \quad }
\alpha\eq 0} \mbox{.} \right.
$$

Now  consider two coordinate functions: the affine space ${\cal
E}_{K}^{(0)}$ has coordinates $a_{ij}$, where $a_{ij}$ are
coefficients of the series $1+\sum_{i,j} a_{ij}u^i\partial^{-j}$
(any element from ${\cal E}_{K}^{(0)}$). For our projective line we
have the following coordinates under the embedding:
$$
a_{-1,1}\eq \left\{ {0\mbox{\quad if\quad} \alpha\ne 0} \atop {-1 \mbox{\quad if\quad}
\alpha\eq 0} \right.
$$
$$
a_{0,1}\eq  \left\{ {-\beta /\alpha \mbox{\quad if\quad} \alpha\ne
0} \atop {0 \mbox{\quad if\quad} \alpha\eq 0}   \mbox{.}  \right.
$$

These functions are not continuous! Hence, it follows that the
embedding of the universal Sato Grassmannian into the affine space
is only a set-theoretic map.

\section{New equations of KP-type on skew fields}

In this section we give the answer on a question set up in the
paper~\cite{P3}. We recall that the classical KP-hierarchy is
constructed by means of the ring of pseudo-differential operators
$P\eq k((x))((\partial^{-1}))$. This ring is a skew field. The point
is to consider other skew fields instead of this one.  We study
below, if there are some new non-trivial generalizations of the
KP-hierarchy for a known list from the classification of
two-dimensional local skew fields. In particular, we give a number
of new equations of the KP-type.

We recall that in the paper~\cite{P3} A.N.Parshin pointed out a
class of non-commutative local fields and showed that these skew
fields possess many properties of commutative fields. He defined a
skew field of formal pseudo-differential operators in $n$ variables
and studied some of their properties. He raised a problem of the
classification of all non-commutative local skew fields. In the
paper~\cite{Zh1} this problem was solved for $n \eq 2$. A.B.~Zheglov
obtained the list of skew fields up to an isomorphism. Among other
skew fields this list contains a classical ring of
pseudo-differential operators $P$ mentioned in the previous part of
this paper. The following theorem was proved in~\cite{Zh1}:

 \begin{th}
\label{main2}
(I) Let $K$ be a two-dimensional local skew field with a commutative residue
skew field.

It splits\footnote{a local skew field is called splittable if there
is a section of the residue homomorphism} if the canonical
automorphism\footnote{the canonical automorphism of a local skew
field is defined as an automorphism of the residue skew field
induced by the inner automorphism $Ad(z)$ of the skew field, where
$z$ is an arbitrary parameter} $\alpha$ satisfies the condition
$\alpha^n\ne Id$ for all $n$. If this condition does not hold, there
are
examples of non-splittable skew fields. \\

(II) Let $K, K'$ be skew fields as in (I). Assume that $\alpha'^n
\ne Id$ for all $n$. Then

(a) $K$ is isomorphic to a two-dimensional local skew field $\bar K
((z))$, where $za \eq a^{\alpha} z$, $a \in \bar K$ and $\bar K$ is
a one-dimensional
local field with the residue field $k$.\\

(b) $K$ and $K'$ are isomorphic if and only if $k\cong k'$, and
there is an isomorphism
$f: \bar K\mapsto \bar K'$ such that $\alpha\eq f^{-1}\alpha' f$.\\

(c) If $char K \eq char k$, $char K' \eq char k'$ and $k,k'$ are
algebraically closed fields of characteristic $0$, then
 $K$ is isomorphic to $K'$ if and only if $k\cong k'$
and $(a_1,i_{\alpha}, y(\alpha ))\eq (a'_1, i_{\alpha'}, y(\alpha' ))$.\\

(III) Let $K, K'$ be splittable two-dimensional local skew fields of
characteristic 0, $k \subset Z(K)$, $k' \subset Z(K')$, and
$\alpha^n \eq Id$, $\alpha'^{n'} \eq Id$ for some natural numbers
$n, n'\ge 1$. Then

(a) $K$ is isomorphic to a two-dimensional local skew field
$k((u))((z))$, where
$$zuz^{-1}\eq
\xi u+u^{\delta'_{i_n}}z^{i_n}+u^{\delta'_{2i_n}}z^{2i_n} \mbox{,}$$
where $\xi^n \eq
1$, $i_n\eq i_n(0,\ldots ,0)$,\\
$\delta'_{i_n}(u) \eq cu^{r_n}$,
$c\in k^*/(k^*)^{e}$, $e \eq (r_n-1,i_n)$,\\
$\delta'_{2i_n}(u) \eq (a_n(0,\ldots
,0)+r_n(i_n+1)/2)u^{-1}(\delta'_{i_n}(u))^2$\\
($i_n, r_n, a_n$ were defined in~\cite{Zh1}).\\
If $n \eq 1$, $i_n \eq \infty$, then $K$ is commutative.\\

(b) $K$ is isomorphic to  $K'$ if and only if  $k\cong k'$ and the sets \\
 $(n,\xi ,i_n,r_n,c,a_n)$,  $(n',\xi' ,i'_n,r'_n,c',a'_n)$
coincide.
\end{th}

\begin{cons}
\label{novoe} Every two-dimensional local skew field $K$ with the
set
$$
(n,\xi ,i_n,r_n,c,a_n)
$$
is a finite-dimensional extension of the skew field with the  set
$(1,1,1,0,1,a)$.
\end{cons}

For every skew field from this list it is possible to define a
decomposition $K\eq K_{+}+K_{-}$, where $K_{-}\eq \{ L\in K:
ord(L)<0\}$\footnote{we use $ord(\sum\limits_k a_k z^k)= - min \{ k
: a_k \ne 0 $ \}}, and $K_{+}$ consists of the operators containing
only nonpositive ($\le 0$) powers of $z$, and a "KP-hierarchy" in
the Lax form:
$$
\frac{\partial L}{\partial t_n}\eq [(L^n)_{+},L] \mbox{,}
$$
where $L\in z^{-1}+K_{-}\otimes k[[\ldots ,t_m,\ldots ]]$. Let $L\eq
z^{-1}+ u_1z+u_2z^2+\ldots $, where $u_m\eq u_m(u, t_1, t_2, \ldots
)$. In the sequel we  will write $\partial /\partial t_n$ as
$\partial_n$.

One can easily  check that if the canonical automorphism $\alpha$
from the classification theorem~\ref{main2} is not trivial, then our
"KP-hierarchy" becomes trivial in a sense that it can be easily
linearized and solvable. We omit calculations here. Thus, we can
assume that $\alpha\eq id$. In this case, if $i>1$, then we have
$[(L^n)_{+},L] \eq -[(L^n)_{-},L] \eq 0 \mbox{\quad mod\quad }
\wp^i$, where $\wp$ is the  ideal of the first valuation in $K$.
Then, as before, our "KP-hierarchy"  is linear and easily solvable.

So, we have to consider the case $i \eq 1$. In this case we have $r
\eq 0$ and $c \eq 1$ (see~\cite{Zh1}), and, therefore, there is only
 one non-trivial parameter $a$. If $a \eq 0$, then $K$ is isomorphic to
the ring $P$ of usual pseudo-differential operators under $x = -u$,
$\partial = z^{-1}$. By $u'_i, u''_i, \ldots $ we denote the
subsequent derivatives by $x$.

For $n\eq 1$ we obtain the following equation
$$
\partial_1u_1\eq u'_1 \mbox{,}
$$
i.e.,  we can assume that $t_1 \eq x$ for $u_1$.

Now we write down the first two equations for $n\eq 2$ and the first equation
for $n\eq 3$.
\begin{equation}
\label{(1)}
\partial_2u_1\eq {u''}_1+2u'_2
\end{equation}
\begin{equation}
\label{(2)}
\partial_2u_2\eq 2u'_3+2u_1u'_1+{u''}_2 - a (2x^{-1}u'_2 +2 x^{-1} u''_1 -x^{-2}u'_1)
\end{equation}
\begin{equation}
\label{(3)}
\partial_3u_1\eq
{u'''}_1+3{u''}_2+3u'_3+6u_1u'_1-3a(x^{-1}{u''}_1-x^{-2}u'_1)
\mbox{.}
\end{equation}
We introduce the following notation: $u\eq u_1(x,y,t)$, where $y \eq
t_2$, $t\eq t_3$. We will use also the standard notations $u_t, u_y,
u_{yy}, \ldots$ for partial derivatives.

Eliminate $u'_3$ from equations (\ref{(2)}) and (\ref{(3)}). We
obtain
\begin{equation}
\label{(4)} 3u_{2y}-2u_t\eq -6uu' -3{u''}_2-2{u'''}- 3a (2x^{-1}u'_2
+ x^{-2}u')  \mbox{.}
\end{equation}
From equation~(\ref{(1)}) we find
$$
{u'''}_2\eq 1/2({u''}_y-u'''') \qquad, \qquad u'_{2y}\eq
1/2(u_{yy}-{u''}_y) \mbox{.}
$$
Differentiating  equation (\ref{(4)}) by $x$ and using these
expressions we finally obtain a new nontrivial KP-equation:
$$
(4u_t-u'''-12uu')'\eq
3u_{yy}+6a(2x^{-2}u''-x^{-2}u_y-x^{-1}u'''+x^{-1}u'_y-2x^{-3}u')
\mbox{.}
$$
It is easily to see that if $a \eq 0$, then we have the usual
KP-equation (compare also explicit calculations in the
paper~\cite{Par}).

Other interesting generalizations of KP-hierarhies which are
connected with the development of Krichever correspondence to higher
dimensions (sections 1-3) are described in~\cite{Zh2}.

\noindent D.V. Osipov \\
Steklov Mathematical Institute RAS\\
e-mail  ${d}_{-} osipov@mi.ras.ru$

\vspace{0.5cm}

\noindent A.B. Zheglov \\
Moscow State Lomonosov University\\
e-mail  $azheglov@mathematik.hu-berlin.de$

\end{document}